\documentclass[12pt,a4paper]{amsart}
\usepackage{amssymb}
\usepackage{mathrsfs}

\headsep=1cm \numberwithin{equation}{section}
\newtheorem{thm}{Theorem}[section]
\newtheorem{cor}[thm]{Corollary}
\newtheorem{lem}[thm]{Lemma}
\newtheorem{prop}[thm]{Proposition}
\theoremstyle{definition}
\newtheorem{defn}[thm]{Definition}
\theoremstyle{remark}
\newtheorem{rem}[thm]{Remark}
\numberwithin{equation}{section}

\newcommand\Ass{\operatorname{Ass}}

\newcommand\Spec{\operatorname{Spec}}
\newcommand\Rad{\operatorname{Rad}}

\begin{document}
\title[Ratliff-Rush closures  and linear growth of primary decompositions]
{Ratliff-Rush closures  and linear growth of primary decompositions of ideals}
\author[Monireh Sedghi]{Monireh Sedghi$^*$}
\address{Department of Mathematics, Azarbaijan Shahid Madani University, Tabriz, Iran. }%
\email{m\_sedghi@tabrizu.ac.ir} \email {sedghi@azaruniv.ac.ir} %
\thanks{ 
This work  has been supported by a  research fund number 217/D/7154 from the Azarbaijan Shahid Madani University. \\
$^*$Corresponding author: e-mail: {\it m\_sedghi@tabrizu.ac.ir} (Monireh Sedghi)}%

\begin{abstract}
Let $R$ be a commutative Noetherian  ring,  $E$  a non-zero finitely 
generated $R$-module and  $I$  an ideal of $R$. One purpose of this paper is to show that the sequences $\Ass_RE/
\widetilde{I_E^n}$ and $\Ass_R\widetilde{I^n _E}/\widetilde{I^{n+1}_E}$,  $n = 1,2, \dots$,  of associated prime ideals
are increasing and eventually stabilize. This extends the main result of
Mirbagheri and Ratliff  \cite[Theorem 3.1]{MR}.  In addition, a characterization
concerning the set $\widetilde{A^*}(I,E)$  is included.  A second purpose of this paper is to prove that $I$ has linear growth primary decompositions for Ratliff-Rush
closures with respect to $E$,  that is, there
exists a positive integer $k$ such that for every positive integer $n$, there exists a
minimal primary decomposition $\widetilde {I^n_E}= Q_1\cap \cdots\cap
Q_s$ in $E$ with $(\Rad(Q_i:_RE))^{nk}\subseteq (Q_i:_R E)$, for all $i= 1, \dots, s$.

\end{abstract}
\maketitle
\section{Introduction}
 Let $I$ be an ideal of a Noetherian ring $R$, and let $E$ be a
non-zero finitely generated module over $R$. We denote by
$\mathscr{R}$ the {\sl Rees ring} $R[u, It] :=\oplus_{n\in
\mathbb{Z}}I^nt^n$ of $R$ with respect to $I$, where $t$ is an indeterminate
and $u =t^{-1}$. Also, the {\sl graded Rees module} $E[u, It] :=
\oplus_{n\in \mathbb{Z}}I^nE$ over $\mathscr{R}$ is denoted by
$\mathscr{E}$, which is a finitely generated graded
$\mathscr{R}$-module. We shall say that $I$ is $E$-{\sl proper} if
$E\neq IE$, and, in this case, we define the $E$-{\sl
grade of} $I$ (written ${\rm grade}\,(I, E)$) to be the maximum length of
all $E$-sequences contained in $I$. In \cite{RR}
Ratliff and Rush studied the interesting ideal,
\begin{center}
$\widetilde{I}= \cup_{n \geqslant1}(I^{n+1}:_R I^n)$= \{$x\in R :
xI^n\subseteq I^{n+1}$ for some $n\geq1$\},
\end{center}
associated with $I$.
 If ${\rm grade}\,I>0$, then this new ideal has some nice properties.
For instance, for all sufficiently large integers $n$, $\widetilde{I^n}=
I^n$. They also proved the interesting fact that, for any
$n\geqslant1$, $\widetilde{I^n}$ is the eventual stable value of the
increasing sequence,

$$(I^{n+1}:_R I)\subseteq (I^{n+2}:_R I^2)\subseteq (I^{n+3}:_R
I^3)\subseteq \cdots.$$

In particular, Mirbagheri and Ratliff, in \cite[Theorem  3.1]{MR} showed that
the sequences of associated prime
ideals \[{\rm Ass}_RR/\widetilde{I^n} \ \ {\rm and}\ \ {\rm
Ass}_R\widetilde{I^n}/\widetilde{I^{n+1}}, \ \ \ n = 1, 2,
\dots\] are increasing and eventually stabilize.
In \cite{HLS}, a regular ideal $I$, i.e., ${\rm grade}\,I>0$, for which
$\widetilde{I} = I$ is called a {\sl Ratliff- Rush ideal}, and the ideal
$\widetilde{I}$ is called the {\sl Ratliff-Rush ideal associated with}
the regular ideal $I$. (For more information about the Ratliff-Rush
ideals, see \cite{HJLS2}, \cite{HLS}, \cite {Na} and \cite{RS}.)
W. Heinzer et al. \cite{HJLS1} introduced a concept analogous to
this for modules over a commutative ring. Let us recall the following definition:

\begin{defn}
(See, Heinzer et al. \cite{HJLS2}). Let $R$ be a commutative ring,
 let $E$ be an $R$-module and let $I$ be an ideal of $R$.
 The Ratliff-Rush closure of $I$ with respect to $E$ denoted by $\widetilde{I_E}$, is
 defined to  be the union of $(I^{n+1}E:_E I^n)$, where $n$ varies in $\mathbb{N}$; i.e.,
 \begin{center}
 $\widetilde{I_E}$=\{$e\in E: I^ne\subseteq I^{n+1}E$ for some $n$\}.
 \end{center}
\end{defn}
If  $E=R$ then the definition reduces to that of the usual Ratliff-Rush ideal
associated to $I$ in $R$ (see \cite{RR}). Furthermore $\widetilde{I_E}$ is a submodule of $E$ and
it is easy to see that $IE\subseteq \widetilde{I}E\subseteq \widetilde{I_E}$. The ideal $I$ is said to be
Ratliff-Rush closed with respect to $E$ if and only if $IE=\widetilde{I_E}$.\\

In the second  section, we study the asymptotic stability of
Ratliff-Rush closures of powers of ideals to modules and, then we extend Mirbagheri and Ratliff's result (see \cite[Theorem  3.1]{MR}). More precisely, we shall show that:
\begin{thm}
Let $R$ be a commutative  Noetherian ring, $I$ an ideal of $R$,  and let $E$ be a non-zero finitely generated $R$-module.
Then the sequences
 \[\Ass_R E/ \widetilde{I^n_E} \ \ {\rm and}\ \ \Ass _R \widetilde{I^n _E}/\widetilde{I^{n+1}_E}
 , \ \ \ n = 1, 2,
\dots\]
 of associated primes are increasing and eventually constant.
\end{thm}

The proof of Theorem 1.2 is given in 2.2 and 2.4.  One our main tools for proving Theorem 1.2 is the following which plays a key role in this section.

\begin{lem}
Let $R$, $I$ and $E$ be as in Theorem {\rm 1.2}. Let $m$ and $n$ be two
natural numbers such that $n\geq m$. Then
$$\widetilde{I^n_E} :_E \widetilde{I^m} = \widetilde{I^n_E}:_EI^m =\widetilde{I^{n-m}_E}.$$
 \end{lem}

Pursuing this point of view further we introduce the notion of Ratliff-Rush reductions of
submodules with respect to an ideal; and it is shown that many of the general properties of standard
reductions extend to Ratliff-Rush reductions. Finally, we will give a characterization of
$\widetilde{A^*}(I, E)$ in terms of the Rees ring and the Rees module of $E$ with respect to $I$.\\

The {\it linear growth} property was first
proved by Irena Swanson \cite{Sw1} and then by R.Y. Sharp using
different methods and in a more general situation in \cite{Sh}. 
We  say  that an ideal $I$ of a Noetherian ring $R$ has  linear  growth  of  primary
decompositions if  there  exists  a  positive  integer $h$ such  that,  for  every
positive integer $n$, there exists a minimal primary decomposition of $I^n$,

\begin{center}
$I^n=\frak q_{n1}\cap \dots \cap \frak q_{nk_n}$, \,\,\,\,\, with $\Rad(\frak q_{ni})^{nh}\subseteq \frak q_{ni}$ for all $i=1, \dots, k_n.$
\end{center}
  I. Swanson  proved  in   \cite[Theorem 3.4]{Sw1}   that  every  proper  ideal  in  every
commutative  Noetherian  ring  has  linear  growth  of  primary  decompositions. Subsequently,  
R.Y. Sharp generalized Swanson's theorem based on the well-known theory of injective $R$-modules and D. Kirby's Artin-Rees Lemma for Artinian modules (see \cite[Proposition 3]{Ki}).

The  purpose  of  the Section 3  is  to  establish  a  parallel  result  for  the Ratliff-Rush closures   of ideals with respect to modules.  We say that  $I$ {\it has linear growth of primary decompositions  for  Ratliff-Rush closures  with respect to $E$} if  there  exists  a  positive  integer $t$ such  that,  for
every integer $n$, there exists a minimal primary decomposition of $\widetilde{I^n_E}$,
\begin{center}
$\tilde {I^n_E}= Q_{n1}\cap \cdots\cap Q_{nl_n}$, \,\,\,\,\, with
 $(\Rad(Q_{ni}:_RE))^{tn}\subseteq (Q_{ni}:_R E)$ for all $i= 1, \dots, l_n$, 
\end{center}

We  prove  in  this  section  that  every  $E$-proper ideal $I$  in  every  commutative
Noetherian  ring  has  linear  growth of primary decompositions for  Ratliff-Rush closures  with respect to  finitely generated $R$-module $E$, whenever ${\rm grade}\,(I, E)>0$. 

Throughout this paper, $R$ will always be a commutative Noetherian
ring with non-zero identity and $I$ will be an ideal of $R$.
For any ideal $J$ of $R$,  the {\it radical} of $J$, denoted by $\Rad(J)$, is defined to
be the set $\{x\in R \,: \, x^n \in J$ for some $n \in \mathbb{N}\}$. 
If $(R, \frak{m})$ is a Noetherian local ring, then $R^\ast$ 
denotes the completion of $R$ with respect to the $\frak{m}$-adic topology.
For any unexplained notation and terminology we refer
 to \cite{N} and \cite{Mat}.
\section{Asymptotic stability of Ratliff-Rush closures of ideals}
One of the main points of this section (see Theorem 2.2) is to generalize and to provide a short proof of the main result of Mirbagheri and Ratliff (see \cite[Theorem  3.1]{MR}); and we will give a characterization of $\widetilde{A^*}(I, E)$ in terms of Rees ring. Further,  we introduce the Ratliff-Rush reductions of
submodules with respect to an ideal; and it is shown that many of the general properties of standard
reductions extend to Ratliff-Rush reductions. The following lemma plays a key role in  this section. 

\begin{lem}
Let $R$ be a commutative Noetherian ring, $I$ an ideal of $R$ and let $E$ be a non-zero finitely generated $R$-module.
Suppose that $m, n$ are two
natural numbers such that $n\geq m$. Then
$$\widetilde{I^n_E} :_E \widetilde{I^m} = \widetilde{I^n_E} :_E I^m =\widetilde{I^{n-m}_E}.$$
 \end{lem}
\proof  Let $e\in  \widetilde{I^n_E} :_E I^m $. Then $$I^me\subset \widetilde{I^n_E}=\bigcup_{k \in \mathbb{N}}
(I^{n+m}E:_EI^k)=I^{n+s}E:_E I^s,$$ for some natural number $s$.
Hence $I^{m+s} e \subset I^{n+s}E$, and so $$ e\in I^{n+s}E:_E I^{m+s}=I^{n-m+m+s}E:_E I^{m+s} \subseteq
\widetilde{I^{n-m}_E}.$$ Now, we show that $\widetilde{I^{n-m}_E}\subseteq \widetilde{I^n_E}:_E \widetilde{I^m}$.
To do this, it is enough for us to prove that $ \widetilde{I^m}\widetilde{I^{n-m}_E}\subseteq \widetilde{I^n_E}$.
Let $e\in \widetilde{I^{n-m}_E}$ and $r\in \widetilde{I^m}$. Then, for some large
$k$, $I^ke\subseteq I^{n-m+k}E$ and $rI^k\subseteq I^{k+m}$.
Hence $$I^k re \subseteq rI^{n-m+k} E=rI^k I^{n-m} E \subseteq I^{k+m}I^{n-m} E =I ^{n+k} E.$$
Therefore $re \in I^{n+k} E:_E I^k \subseteq \widetilde{I^n_E}$, as required.  \qed \\

Now we are prepared to prove the main result of this section, which
is an extension of Mirbagheri and Ratliff's result.
\begin{thm}
Let $R$ be a commutative Noetherian ring and let $E$ be a non-zero finitely generated $R$-module. Suppose
that $I$ is an ideal of $R$.
 Then the sequence
 $\{\Ass_R E/ \widetilde{I^n_E}\}_ {n \in \mathbb{N}}$, of associated primes, is increasing and eventually constant.
\end{thm}
\proof Let $\frak{p}\in \Spec R$ and $\frak{p}=\widetilde{I^n_E} :_R e$, for some $e\in E$. Then, in view of the Lemma 2.1,
$\widetilde{I^n_E} =\widetilde{I^{n+1}_E}:_R I$, and so $\frak{p}= \widetilde{I^{n+1}_E}:_R I e$.
Since $I$ is finitely generated, we have $\frak{p}= (\widetilde{I^{n+1}_E}:_Rf)$ for some $f \in Ie$.
Therefore $\frak{p} \in Ass_R E/ \widetilde{I^{n+1}_E}$. This shows the sequence $\{\Ass_R E/ \widetilde{I^n_E}\}_ {n \in \mathbb{N}}$
is increasing.

On the other hand,  $\widetilde{I^n_E}= (I^{n+s}E:_E I^s)$ for some $s\in\mathbb{N}$.
Then we have $$\frak{p}= (I^{n+s} E:_E I^s) :_R e =I^{n+s}E:_ R I^s e,$$
and hence $\frak{p} \in \Ass_R E/I^{n+s}E $. Consequently,
$$ \bigcup _{n \geq 1} \Ass_R E/ \widetilde{I^n _E } \subseteq \bigcup _{k \geq 1} \Ass_R E/I^kE.$$
Now the desired result follows from Brodmann's Theorem (cf. \cite{B}).   \qed \\

\begin{defn}
Let $R$ be a commutative Noetherian ring, $E$  a non-zero finitely generated $R$-module,
and let $I$ be an ideal of $R$. Then the eventual constant value of the sequence
$\Ass_R E/ \widetilde{I^n_E}$, $n=1, 2, \dots,$ will be denoted by
$\widetilde{A^*}(I, E)$.
\end{defn}

\begin{prop}
Let $R$ be a commutative Noetherian ring, $E$ a non-zero finitely generated $R$-module,
and let $I$ be an ideal of $R$.  Then the sequence $\{\Ass_R\widetilde{I^n _E}/\widetilde{I^{n+1}_E}\}_{n \geq 1}$
is monotonically increasing and eventually stable.
\end{prop}
\proof  Let $\frak{p} \in \Ass_R \widetilde{I^n _E}/\widetilde{I^{n+1}_E}$. Then there is $e\in \widetilde{I^n _E}$  such that,
$$\frak{p}= \widetilde{I^{n+1}_E}:_Re =(\widetilde{I^{n+2}_E}:_E I):_R e = \widetilde{I^{n+2}_E}:_R Ie.$$
Since $Ie\subseteq \widetilde{I^{n+1}_E}$, it yields that $\frak{p} \in \Ass_R \widetilde{I^{n+1}_E}/ \widetilde{I^{n+2}_E}$.
Now we can process similarly to the proof of Theorem 2.2 to deduce that the set
$\Ass_R \widetilde{I^n _E}/\widetilde{I^{n+1}_E}$ is increasing and eventually constant.  \qed \\

\begin{cor}
{\rm(cf. \cite[Theorem  3.1]{MR})} Let $R$ be a commutative Noetherian ring
and $I$ an ideal of $R$. Then,  for
all  $n\in \mathbb{N}$,  $${\rm Ass}_RR/\widetilde{I^{n}} = {\rm Ass}_R\widetilde{I^{n-1}}/\widetilde{I^{n}},$$  
and   the sequences of associated prime
ideals \[{\rm Ass}_RR/\widetilde{I^n} \ \ {\rm and}\ \ {\rm
Ass}_R\widetilde{I^n}/\widetilde{I^{n+1}}, \ \ \ n = 1, 2, \dots\]  are increasing and eventually stabilize.
\end{cor}
\proof Clearly ${\rm Ass}_R\widetilde{I^{n-1}}/\widetilde{I^{n}} \subseteq {\rm Ass}_RR/\widetilde{I^{n}}.$ Let $\frak p\in {\rm Ass}_RR/\widetilde{I^{n}}$. Then there exists $r\in R$ such that $\frak p=(\widetilde{I^{n}}:_Rr)$.  Since $I\subseteq \widetilde{I}\subseteq \frak p$,
it follows that $Ir\subseteq \widetilde{I^{n}}$, and so $r\in (\widetilde{I^{n}}:_R I)$.  Hence, in view of Lemma 2.1,  $r \in  \widetilde{I^{n-1}}$. Thus 
$\frak p\in {\rm Ass}_R\widetilde{I^{n-1}}/\widetilde{I^{n}}$ and ${\rm Ass}_RR/\widetilde{I^{n}} = {\rm Ass}_R\widetilde{I^{n-1}}/\widetilde{I^{n}}$. Now, the result follows from Theorem 2.2. \qed \\

In the next result, we introduce the concept of a Ratliff-Rush reduction of submodules of an
$R$-module $E$ with respect to an ideal $I$ of $R$, and then show that most of the known basic properties of standard reductions of ideals (see \cite{NR}) also
hold for Ratliff-Rush reductions.\\
\begin{defn}
Let $R$ be a commutative ring and $E$ an $R$-module.  Let $N_1\subseteq N_2$ be submodules of $E$,
and let $I$ be an ideal of $R$. We say that $N_1$ is
a {\it Ratliff-Rush reduction of} $N_2$ with respect to $I$ if $ N_2 \subseteq \widetilde{I_{N_1}}$, (note that $\widetilde{I_{N_1}}\subseteq \widetilde{I_{N_2}}\subseteq \widetilde{I_{E}}$). In other word, for any $x\in N_2$
 there exists a  $n\in \mathbb{N}$ such that $I^nx \subseteq I^{n+1} N_1$.
\end{defn}

\begin{rem}
Let $R$ be a commutative Noetherian ring and $E$ a finitely generated $R$-module.
Let $N_1\subseteq N_2$ be submodules of $E$ such that $N_1$ is a Ratliff-Rush reduction of $N_2$ with respect to $I$.
Then there exists a $n\in \mathbb{N}$ such that  $I^nN_2 \subseteq I^{n+1}N_1$.
\end{rem}

\begin{thm}
Let $R$ be a commutative ring and $E$ an $R$-module. Let $I$ be an ideal of $R$ and let $N,K$ and $L$
be submodules of $E$. Then the following conditions hold:

{\rm(i)} If  $R$ is Noetherian  and $E$ is finitely generated such that $N\subseteq K\subseteq L$
and $N$ is a reduction of $K$ with respect to $I$ and $K$ is a reduction of $L$ with respect to $I$; then $N$ is a reduction of $L$ with respect to $I$.

{\rm(ii)} If $N\subseteq K $ then $\widetilde{I_N}\subset\widetilde{I_K}$.

{\rm(iii)} If $N\subseteq K$, $N$ is a reduction of $K$ with respect to $I$ and  $J$ is an ideal of $R$, then
 $JN$ is a reduction of  $JK$ with respect to $I$.

{\rm(iv)} If $N_1\subseteq N_2 \subseteq N_3$ are submodules of $E$ such that $N_1$ is a reduction of $N_3$ with respect to $I$,
then  $N_2$  is  a reduction of $N_3$ with respect to $I$.

{\rm(v)} If $N_1$, $N_2$, $M_1$ and $ M_2$ are submodules of $E$ such that $N_1$ is a reduction of $N_2$ with respect to $I$ and  $M_1$
is a reduction of $M_2 $ with respect to $I$,  then $N_1+M_1$ is a reduction of $N_2+M_2$ with respect to $I$.

{\rm(vi)} If $c\in R$ is an $E$-regular element and $N\subseteq K$ such that $cN$ is a reduction of $cK$ with respect to $I$,
then $N$ is a reduction of $K$ with respect to $I$.

{\rm (vii)} If $S$ is a multiplicatively closed subset of $R$ s.t. $0\notin S$, $N$ is a reduction of $K$ with respect to $I$,
then $S^{-1}N$ is a reduction of $S^{-1}K$ with respect to $S^{-1}I$

{\rm(viii)} Let $(R,\frak{m})$ be local (Noetherian) and let $E$ be finitely generated. The $K \subseteq L$ is
 a reduction of $L$ with respect to $I$ if and only if $KR^*$ is a reduction of $LR^*$ with respect to $IR^*$.
\end{thm}

\proof  (i) Let $x\in L$ be an arbitrary element. Then since $K$ is a reduction of $L$ with respect to $I$,
then $I^n x\subseteq I^{n+1}K$ for some $n\in \mathbb{N}$. Also, since $N$ is a reduction of $K$ with respect to $I$ and $R$
and $E$ are Noetherian, it follows that there exists a $m \in \mathbb{N}$, such that $I^mK \subseteq I^{m+1}N$.
Hence $N$ is a reduction of $L$ with respect to $I$.

(ii) Let $x \in \widetilde{I_N}$. Then $I^n x\subseteq I^{n+1} N$ for some $n\in \mathbb{N}$.
Since $N\subseteq K$, it yields that $I^{n+1} N \subseteq I^{n+1}K$. Thus $I^n x\subseteq I^{n+1}K$, and so $x\in \widetilde{I_K}$.

(iii) Let  $N$ be a reduction of $K$ and $J$ be an ideal of $R$. Then, for any $x\in K$, there
is $n\in \mathbb{N}$ such that $I^nx \subseteq I^{n+1}K$. Hence  $I^n Jx \subseteq I^{n+1}JN$, and so
$JN$ is a reduction of $JK$ with respect to $I$.

(iv)  Let $y \in N_3$ be an arbitrary element. As $N_1$ is a reduction of $N_3$ with respect to $I$, there
 exists $n\in \mathbb{N}$ such that $I^ny\subseteq I^{n+1}N_1$. Now, since $N_1 \subseteq N_2$,
it follows that $I^ny \subseteq I^{n+1}N_2$. Hence $N_2$ is a reduction of $N_3$ with respect to $I$.

(v) Let $x\in N_2$ and $y\in M_2$ be arbitrary elements. Then there is $n\in\mathbb{ N}$ such that
$I^n x \subseteq I^{n+1}N_1$ and $I^n y \subseteq I^{n+1}M_1$. Whence $I^n (x+y)\subseteq I^{n+1}(N_1 +M_1)$,
and so $N_1+M_1$ is a reduction of $N_2+M_2$ with respect to $I$.

(vi)  Let $x \in N_2$. Then  $cx \in \widetilde{I_{cN_1}}$, and so there is $n\in\mathbb{ N}$ such that
$I^n cx \subseteq I^{n+1} cN_1$.  Now, let $a \in I^n$ be an arbitrary element. Then $acx=cz$ for some $z \in I^{n+1}N_1$,
and so $N_1$  is a reduction of $N_2$ with respect to $I$.

(vii) Let $x \in N_2$. Then $I^nx \subseteq I^{n+1}N_1$ for some  natural number $n$. Hence
 $$(S^{-1}I)^n x \subseteq (S^{-1}I)^{n+1} (S^{-1}N_1),$$ and so $S^{-1}N_1$ is a reduction of  $S^{-1}N_2$ with respect to $S^{-1}I$.

(ix) Let $x \in L$. Then $I^nx \subseteq I^{n+1}K$ for some  natural number $n$. Thus
 $$(IR^*)^nx \subseteq (IR^*)^{n+1}(KR^*),$$ and so $KR^*$ is a reduction of $LR^*$ with respect to $IR^*$.
Since $R^*$ is a faithfully flat $R$-module the converse holds too.  \qed \\

We end this section with a characterization of $\widetilde{A^*}(I,E)$ in terms of Rees ring.
\begin{thm}
 Let $R$ be a commutative Noetherian ring and let $E$ be a non-zero finitely generated $R$-module. Suppose that  $I$ is
an $E$-proper ideal of $R$ such that ${\rm grade} \, (I, E)>0$. Then the following statements are equivalent:

{\rm(i)} $\frak{p} \in \widetilde{A^*}(I, E)$.

{\rm(ii)} There exists a prime ideal $\frak q \in \widetilde{A^*}(t^{-1}\mathscr{R}, \mathscr{E})$
such that $\frak q \cap R= \frak{p}$.
\end{thm}

\proof $\rm(i)\Longrightarrow(ii)$. Let  $\frak{p} \in \widetilde{A^*}(I, E)$. Then
there exists an integer $n\geq1$ such that $\frak p\in \Ass_R E/ \widetilde{I^n_E}$. Now,
in view of \cite[Lemma 2.1]{Na}, we have  $$\widetilde{I^n_E}= I^{n+r}E:_E I^r= I^nE,$$
for some large integer $n$ and for every integer $r\geq1$. Hence $\frak p\in \Ass_R E/ I^nE$. Since $I^nE= t^{-n}\mathscr{E}\cap E$,
it follows that there is a prime ideal  $\frak q \in \Ass_{\mathscr{R}} \mathscr{E}/t^{-n}\mathscr{E}$
such that $\frak q \cap R= \frak{p}$. Now, as
$$(\widetilde{t^{-n}\mathscr{R})_{\mathscr{E}}}=t^{-(n+r)}\mathscr{E}:_{\mathscr{E}} t^{-r}\mathscr{R}$$
it is easy to see that $\frak q \in \Ass_{\mathscr{R}} \mathscr{E}/(\widetilde{t^{-n}\mathscr{R}})_{\mathscr{E}}$.
Therefore $\frak q \in \widetilde{A^*}(t^{-1}\mathscr{R}, \mathscr{E})$
such that $\frak q \cap R= \frak{p}$, as required.

In order to prove the implication $\rm(ii) \Longrightarrow (i)$, suppose $\frak q$ satisfies in (ii).
Then by definition $\frak q\in \Ass_{\mathscr{R}} \mathscr{E}/ (\widetilde{t^{-n}\mathscr{R})_{\mathscr{E}}}$  for large $n$.
Now, in view of \cite[Lemma 2.1]{Na}, we have
 $$(\widetilde{t^{-n}\mathscr{R})_{\mathscr{E}}}= t^{-(n+r)}\mathscr{E}:_{\mathscr{E}} t^{-r}\mathscr{R}= t^{-n}\mathscr{E},$$
for some large integer $n$ and for every integer $r\geq1$. Whence $\frak q\in \Ass_{\mathscr{R}}\mathscr{E}/t^{-n}\mathscr{E}$. Now, 
it is easy to see that if $Q$ is a $\frak q$-primary component of $t^{-n}\mathscr{E}$ in $\mathscr{E}$, then $Q\cap E$ is a  $\frak p$-primary component of 
$I^nE$ in $E$. Therefore, as 
$$(\widetilde{t^{-n}\mathscr{R})_{\mathscr{E}}}=t^{-(n+r)}\mathscr{E}:_{\mathscr{E}} t^{-r}\mathscr{R}= t^{-n}\mathscr{E},$$ 
 it  follows from $I^nE= t^{-n}\mathscr{E}\cap E$ that $\frak p\in \Ass_RE/I^nE$. Whence  by using again \cite[Lemma 2.1]{Na}, we have $\frak p\in \Ass_R E/ \widetilde{I^n_E}$.
Thus $\frak{p}\in \widetilde{A^*}(I, E)$, as required. \qed\\

\section{ linear growth of primary decompositions}

As the main result we shall prove that, for certain $E$-proper ideal $I$ of $R$,
$I$ has {\it linear growth primary decompositions for Ratliff-Rush closures} with respect to $E$. We begin with:
 \begin{rem}

   Let $A$ be a commutative Noetherian ring and $X$ a finitely
   generated $A$-module. For an ideal $J$ of $A$ and a submodule
   $Y\subseteq X$ the increasing sequence of submodules,

   $$Y\subseteq (Y:_X J)\subseteq (Y:_X J^2)\subseteq \cdots \subseteq
(Y:_X
   J^n) \subseteq \cdots,$$
   becomes stationary. Denote its ultimate constant value by $Y:_X
   \langle J\rangle$. Note that $Y:_X \langle J\rangle = Y:_X J^n$
   for all large $n$.

   One has $Ass_AX/(Y:_X \langle J\rangle)= Ass_A(X/Y) \diagdown
   V(J)$. Therefore the primary decomposition of $Y:_X \langle
J\rangle$
   consists of those primary components of $Y$ whose associated prime
   ideals do not contain $J$.
\end{rem}
We are now ready to state and prove the main theorem of this
section, which shows that an ideal $I$ in a Noetherian ring $R$ has  linear growth primary decompositions for Ratliff-Rush closures with respect to a finitely generated module over $R$.

\begin{thm}
Let $R$ be a Noetherian ring and let $E$ be a non-zero finitely
generated $R$-module. Suppose that $I$ is an $E$-proper ideal of $R$
such that $grade \,(I, E) >0$. Then the following
conditions hold:

{\rm(i)} For any integer $n\geq1$, $(u^n\mathscr{E}:_\mathscr{E}
\langle It\rangle) \cap  E= \widetilde{I^n_E}$.

{\rm(ii)} $I$ has linear growth of primary decompositions for
Ratliff-Rush closures with respect to $E$.

More precisely, there exists a positive integer $t$ such that, for
every $n\in\mathbb{N}$, there exists a minimal primary decomposition
of $\widetilde{I^n_E}$,

$$\widetilde{I^n_E}= Q_{(n,1)} \cap Q_{(n,2)} \cap \cdots \cap
Q_{(n,l_n)},$$

such that $(\Rad(Q_{(n,i)}:_R E))^{tn} \subseteq (Q_{(n,i)}:_R E)$
for all $i= 1, 2, ..., l_n$.

\end{thm}{\it Proof}. (i) Let $n\in\mathbb{N}$. Then, by Remark
3.1, $$(u^n\mathscr{E}:_\mathscr{E} \langle It \rangle)=
(u^n\mathscr{E}:_\mathscr{E} I^kt^k)$$ for all large $k$. On the
other hand, it is easy to see that $$(u^n\mathscr{E}:_\mathscr{E}
I^kt^k)\cap E= (I^{n+k}E:_E I^k).$$ Now, the claim follows from
\cite[Proposition 2.2]{Na}.

(ii) In view of \cite[Theorem 2.1]{Sh}, there exists a positive
integer $t$ such that, for every $n\in\mathbb{N}$, there exists a
minimal primary decomposition of $u^n\mathscr{E}$ in $\mathscr{E}$,
$$u^n\mathscr{E}= Q_{(n,1)} \cap  Q_{(n,2)}\cap \cdots \cap
Q_{(n,l_n)},$$ such that $$(\Rad(Q_{(n,i)}:_\mathscr{R}
\mathscr{E}))^{tn} \subseteq (Q_{(n,i)}:_\mathscr{R} \mathscr{E})$$
for all $i= 1, 2, \dots, l_n$. Now, for $i= 1, 2, \dots, l_n$, we let
$\Rad(Q_{(n,i)}:_\mathscr{R} \mathscr{E})= \frak{p_i}$. After an
appropriate reordering of the $\frak{p_i}$'s, there will be an
integer $s\geq 1$ such that $s\leq l_n$ and $It\subseteq \frak{p_i}$
for $i= 1, 2, \dots, s$ and $It\nsubseteq  \frak{p_j}$ for $j= s+1,
..., l_n$. Then, by Remark 3.1, $$Q_{(n,s+1)} \cap Q_{(n,s+2)}\cap
\cdots \cap Q_{(n,l_n)}$$ is a primary decomposition of
$(u^n\mathscr{E}:_\mathscr{E} \langle It \rangle)$ in $\mathscr{E}$.
Consequently, by virtue of (i), we have $$\widetilde{I^n_E}= (Q_{(n,s+1)}
\cap E) \cap  (Q_{(n,s+2)} \cap E)\cap \cdots \cap (Q_{(n,l_n)} \cap
E).$$ Now the desired result follows easily from this. $\Box$

\begin{center}
{\bf Acknowledgments}
\end{center}
The author is deeply grateful to the referee for a very careful
reading of the manuscript and many valuable suggestions in improving
the quality of the paper. Also, I am  grateful to the Professor Reza Naghipour
 for his careful reading of the first draft. Finally, the author would like to thank the Azarbaijan Shahid Madani University for the
financial support.


\end{document}